\documentclass{article}
\usepackage{amssymb}
\usepackage{amsfonts}
\usepackage{amsmath}
\usepackage{graphicx}

\newtheorem{claim}{Claim}

\newenvironment{definition}[1][Definition]{\medskip \noindent\textbf{#1.} } {\medskip}

\newenvironment{proof}[1][Proof]{\noindent\textbf{#1.} }{\ \rule{0.5em}{0.5em} \medskip}
\def\QATOP#1#2{{#1 \atop #2}}
\def\dsum{\mathop{\displaystyle \sum }}
\def\QOVERD#1#2#3#4{{#3 \overwithdelims#1#2 #4}}
\def\dint{\mathop{\displaystyle \int}}

\begin{document}

\title{Longest common subsequences and the Bernoulli matching model:
numerical work and analyses of the r-reach simplification}
\author{Jonah Blasiak}
\maketitle

\begin{abstract}
The expected length of longest common subsequences is a problem that has
been in the literature for at least twenty five years. Determining the
limiting constants $\gamma_k$ appears to be quite difficult, and the
current best bounds leave much room for improvement. Boutet de Monvel
explores an independent version of the problem he calls the Bernoulli
Matching model. He explores this problem and its relation to the longest
common subsequence problem. This paper continues this pursuit by focusing on
a simplification we term r-reach. For the string model, $\mathbf{L}_{r}(u,v)$
is the longest common subsequence of $u$ and $v$ given that each matched
pair of letters is no more than r letters apart.
\end{abstract}

\section{Introduction}

In our technology oriented society fast processing of digital data is
becoming increasingly important. String comparison is a kind of data
processing that has applications in a wide range of fields including
molecular biology, human speech recognition, computer spelling correction,
and gas chromatography \cite{SK}. A robust, extensively studied, method for
comparing two strings, $u$ and $v$ say, is to compute the length of one of
their longest common subsequences (denote this length by $\mathbf{L}(u,v)$).
A subsequence of a string $u$ is a string obtained by deleting some elements
of $u$. For example, \emph{netra} is a subsequence of \emph{cinematography}.
A longest common subsequence of two strings $u$ and $v$ is a subsequence of $%
u$ and $v$ of maximum length. For example, \emph{netra} is an longest common
subsequence of \emph{cinematography} and \emph{neurotransmitter} because
there is no longer string that is a subsequence of both strings.

\subsection{The Random String model}

The following notation will be useful for working with strings:

\begin{definition}
Define an alphabet $\Sigma $ of size $k$ to be $\{0,1,\ldots ,k-1\}$. Let $%
\Sigma ^{n}$ be the set of all sequences of length $n$ on alphabet $\Sigma $.
\end{definition}

\begin{definition}
If $u=u_{1}u_{2}\ldots u_{n}$ and $u_{i}$ $\in $ $\Sigma $, define $%
u(i\ldots j)$ to be the substring $u_{i}u_{i+1}\ldots u_{j\text{.}}$
\end{definition}

A very interesting and difficult problem is to compute the average length of
longest common subsequences over all possible pairs of \ strings. Or more
precisely, define%
\begin{equation*}
\mathbf{EL}_{n}^{(k)}={\frac{1}{k^{2n}}}\sum_{u,v\in \Sigma ^{n}}\mathbf{L}%
(u,v)
\end{equation*}%
An open problem is to compute the following limit:%
\begin{equation*}
\gamma _{k}=\lim_{n\rightarrow \infty }\frac{\mathbf{EL}_{n}^{(k)}}{n}
\end{equation*}%
Klarner and Rivest established that $\mathbf{EL}_{n}$ is superadditive--$%
\mathbf{EL}_{n+m}\geq \mathbf{EL}_{n}+\mathbf{EL}_{m}$--and from this it can
be shown that the above limit exists (see e.g., \cite{BGNS}).

The current best lower and upper bounds as well as Monte Carlo
approximations of $\gamma _{k}$ are shown in Figure (\ref{current best}).

\begin{figure}[tbp]
\begin{tabular}{l|lll||l|lll}
k & $\QATOP{\text{lower}}{\text{bound}}$ & {\small approximation} & $\QATOP{%
\text{upper}}{\text{bound}}$ & k & $\QATOP{\text{upper}}{\text{bound}}$ &
{\small approximation} & $\QATOP{\text{lower}}{\text{bound}}$ \\ \hline
2 & .77391 & .8123 & .83763 & 9 & .40321 & .4936 & .55394 \\
3 & .63376 & .7176 & .76581 & 10 & .38656 & .4747 & .53486 \\
4 & .55282 & .6544 & .70824 & 11 & .37196 & .4580 & .51785 \\
5 & .50952 & .6075 & .66443 & 12 & .35899 & .4432 & .50260 \\
6 & .47169 & .5707 & .62932 & 13 & .34737 & .4297 & .48880 \\
7 & .44502 & .5405 & .60019 & 14 & .33687 & .4176 & .47620 \\
8 & .42237 & .5152 & .57541 & 15 & .32732 & .4066 & .46462%
\end{tabular}%
\caption{Current best bounds and Monte Carlo approximations of $\protect%
\gamma _{k}$. Lower bounds are from \protect\cite{D} and \protect\cite{BGNS}%
. Upper bounds are from \protect\cite{D}. Approximations are from
\protect\cite{B} and were computed using Monte Carlo simulations
extrapolated to large $n$ using $\frac{\mathbf{EL}n}{n}=\protect\gamma _{k}+%
\frac{A_{k}}{\protect\sqrt{n}\ln n}+\frac{C_{k}}{n\ln n}$, for real numbers $%
A_{k}$, $C_{k}$ that don't depend on $n$.\ }
\label{current best}
\end{figure}

Longest common subsequence computations can also be formulated as a dynamic
programming algorithm or as a directed time passage percolation model (see
e.g. \cite{D},\cite{B}). In the directed time passage percolation model, we
work with the two dimensional lattice in the first quadrant: vertices exist
at the points $(i,j)$ for $i,j\in \{0,1,2,...\}$. On each vertex $(i,j)$ $%
\mathbf{D}_{i,j}$ will is an integer, and $\mathbf{D}_{i,0}$ and $\mathbf{D}%
_{0,i}$ are initialized to $0$. Given two strings $u$ and $v$, $\mathbf{L}%
(u,v)$ is computed by preserving $\mathbf{D}_{i,j}=\mathbf{L}(u(1\ldots
i),v(1\ldots j))$. The necessary recurrence is%
\begin{equation*}
\mathbf{D}_{i,j}=\left\{
\begin{array}{ll}
\mathbf{D}_{i-1,j-1}+1 & \text{if }\delta _{u(i),v(j)}=1 \\
max\{\mathbf{D}_{i,j-1},\mathbf{D}_{i-1,j}\} & \text{if }\delta
_{u(i),v(j)}=0%
\end{array}%
\right.
\end{equation*}%
Where $\delta _{u(i),v(j)}$ is the Kronecker delta (the motivations for this
notation will become clear in the next section). Another way of looking at
this recurrence is to make bonds between adjacent vertices in the lattice
directed in the positive $x$ and $y$ directions. A diagonal bond from $%
(i-1,j-1)$ to $(i,j)$ is added if and only if $\delta _{u(i),v(j)}=1$. If
the horizontal and vertical bonds are given weight $0$, and the diagonal
bonds are given weight $1$, $\mathbf{L}(u,v)$ is the weight of a maximum
weight path from $(0,0)$ to $(|u|,|v|)$.

\subsection{The Bernoulli Matching model}

A related problem called the Bernoulli Matching model is named and well
explored by Boutet de Monvel in \cite{B}. It is most readily seen as a
modification of the directed time passage percolation model. Instead of
placing diagonal bonds based on a match in a pair of strings, diagonal bonds
are placed independently at each location with probability $1/k$. In the
random string model, the probability of a bond between $(i-1,j-1)$ and $%
(i,j) $ is $1/k$, but these probabilities are not independent. The
recurrence for the Bernoulli Matching model is
\begin{equation*}
\mathbf{D}_{i,j}=\left\{
\begin{array}{ll}
\mathbf{D}_{i-1,j-1}+1 & \text{if }\epsilon _{ij}=1 \\
max\{\mathbf{D}_{i,j-1},\mathbf{D}_{i-1,j}\} & \text{if }\epsilon _{ij}=0%
\end{array}%
\right.
\end{equation*}

where the $\epsilon _{ij}$ are independent random variables with $\Pr
(\epsilon _{ij}=1)=1/k$ and $\Pr (\epsilon _{ij}=0)=1-1/k$. Let $\mathbf{EL}%
_{n}^{B(k)}$ be the expected value of $\mathbf{D}_{n,n}$ given this model. $%
\mathbf{EL}_{n}^{B(k)}$, like $\mathbf{EL}_{n}^{(k)}$, is superadditive \cite%
{B} and therefore the following limit exists:
\begin{equation*}
\gamma _{k}^{B}=\lim_{n\rightarrow \infty }\frac{\mathbf{EL}_{n}^{B(k)}}{n}
\end{equation*}

Boutet de Monvel \cite{B} has conjectured that $\gamma _{k}^{B}=\frac{2}{1+%
\sqrt{k}}$ and gives a more general conjecture for the off diagonal lattice
positions (Steele conjectured this for the Random String model in 1982,
Boutet de Monvel refined it in 1999). He also presents a nice derivation of
this result based on cavity methods typically used for the mean field theory
of disordered systems, which he does not try to justify rigorously. Though
not yet a proof, the method appears to solve the problem quite elegantly and
agrees well with numerical approximations.

\subsection{The r-reach simplification}

A straight-forward way of obtaining a lower bound for $\mathbf{EL}_{n}^{(k)}$
is to only consider common subsequences that do not match letters "too far"
from each other. This is equivalent to restricting the lattice to a diagonal
band of fixed width with center line $x=y$. More precisely, let $\mathbf{L}%
_{r}(u,v)$ be the length of a common subsequence of $u$ and $v$ as long as
possible given that if $u(i)=v(j)$ are paired by the subsequence, then $%
|i-j|\leq r$. We will use $\mathbf{R}$ instead of $\mathbf{D}$ when working
with r-reach. The recurrence is modified as follows ($\mathbf{R}_{i,0}$, and
$\mathbf{R}_{0,i}$ are initialized to $0$ as before):%
\begin{equation*}
\mathbf{R}_{i,j}=\left\{
\begin{array}{lll}
\mathbf{R}_{i-1,j-1}+1 & \text{if }\delta _{u(i),v(j)},\left( \epsilon
_{ij}\right) =1 &  \\
\max \{\mathbf{R}_{i,j-1},\mathbf{R}_{i-1,j}\} & \text{if }\delta
_{u(i),v(j)},\left( \epsilon _{ij}\right) =0 & \text{and }|i-j|<r \\
\mathbf{R}_{i,j-1} & \text{if }\delta _{u(i),v(j)},\left( \epsilon
_{ij}\right) =0 & \text{and }j-i\geq r \\
\mathbf{R}_{i-1,j} & \text{if }\delta _{u(i),v(j)},\left( \epsilon
_{ij}\right) =0 & \text{and }i-j\geq r%
\end{array}%
\right.
\end{equation*}

Let $\mathbf{EL}_{n,k,r}$, $\left( \mathbf{EL}_{n,k,r}^{B}\right) $ be the
expected value of $\mathbf{R}_{n,n}$ given this model. Superadditivity still
holds in this model $($i.e. $\mathbf{EL}_{n,k,r}+\mathbf{EL}_{m,k,r}\leq
\mathbf{EL}_{(n+m),k,r})$ because a maximum weight path from $(0,0)$ to $%
\mathbf{(}n+m,n+m)$ has weight at least as large as (weight of maximum
weight path from $(0,0)$ to $(n,n)$)+(weight of maximum weight path from $%
(n,n)$ to $(n+m,n+m)$). The same argument applies to $\mathbf{EL}%
_{n,k,r}^{B} $. Now define%
\begin{equation*}
\gamma _{k,r}=\lim_{n\longrightarrow \infty }\frac{\mathbf{EL}_{n,k,r}}{n}%
\text{, }\gamma _{k,r}^{B}=\lim_{n\longrightarrow \infty }\frac{\mathbf{EL}%
_{n,k,r}^{B}}{n}
\end{equation*}

A simple but quite interesting fact is

\begin{claim}
\begin{equation*}
\lim_{r\longrightarrow \infty }\gamma _{k,r}^{B}=\gamma _{k}^{B}\text{ and }%
\lim_{r\longrightarrow \infty }\gamma _{k,r}=\gamma _{k}
\end{equation*}
\end{claim}

\begin{proof}
r-reach effectively reduces the allowable paths. It is easy to see that for
fixed values of $\epsilon _{ij}$, $\mathbf{D}_{n,n}\geqslant \mathbf{R}%
_{n,n} $, and therefore
\begin{equation*}
\mathbf{EL}_{n,k,r}^{B}\leq \mathbf{EL}_{n}^{B(k)}\Longrightarrow \gamma
_{k,r}^{B}\leq \gamma _{k}^{B}\Longrightarrow \lim_{r\longrightarrow \infty
}\gamma _{k,r}^{B}\leq \gamma _{k}^{B}
\end{equation*}%
Next apply superadditivity and $\mathbf{EL}_{r,k,r}^{B}=\mathbf{EL}%
_{r}^{B(k)}$ to show
\begin{equation*}
\gamma _{k,r}^{B}=\lim_{n\longrightarrow \infty }\frac{\mathbf{EL}%
_{n,k,r}^{B}}{n}\geqslant \frac{\mathbf{EL}_{r,k,r}^{B}}{r}=\frac{\mathbf{EL}%
_{r}^{B(k)}}{r}.
\end{equation*}%
Taking the limit of both sides yields
\begin{equation*}
\lim_{r\longrightarrow \infty }\gamma _{k,r}^{B}\geqslant
\lim_{r\longrightarrow \infty }\frac{\mathbf{EL}_{r}^{B(k)}}{r}=\gamma
_{k}^{B}
\end{equation*}%
The analogous result for the Random String model is proved the same way.
\end{proof}

\section{Solutions to Bernoulli Matching model r-reach for small r}

For small $r$, the percolation problem can be dissected in full detail. The
approach used is fairly straight-foward and computationally intensive.
Unfortunately it appears that the r-reach problem is not as elegant as the
original--possibly because of the "discontinuous" boundary effects at the
displaced diagonals $(i,i+r)$ and $(i+r,i)$. There are several reasons this
problem is worth studying, however. First of all it gives lower bounds for
the original problem. Also, it is an interesting setting to compare the
Random String model with the Bernoulli Matching model. The methods outlined
below seem very difficult to use to solve the problem for general $r$,
however they provide foundations for numerical work on large $r$.

The basic idea of the following analyses is to break the lattice into
sections consisting of the $2r+1$ vertices $%
(n-r,n),(n-r+1,n),...(n,n),(n,n-1),...(n,n-r)$ and then compute
probabilities that $\mathbf{R}$ takes on specific values at these vertices.
We only need to know the distribution of the $n^{th}$ section to compute the
distribution of the $(n+1)^{st}$ section. More formally, let $P_{n}(z)$ be
the probability that $\mathbf{R}_{n,n}=z$. For notational convenience let $%
x_{0}=y_{0}=z$. For $(n\geqslant r)$ let $%
R_{n}(z,x_{1},y_{1},x_{2},y_{2},...,x_{r},y_{r})$ be the event that $(%
\mathbf{R}_{n-i,n}=x_{i\text{ }}$and $\mathbf{R}_{n,n-i}=y_{i}$ $\forall
i\in \{0,1,...,r\})$. Also define
\begin{equation*}
P_{n}(z,x_{1},y_{1},x_{2},y_{2},...,x_{r},y_{r})=\Pr
(R_{n}(z,x_{1},y_{1},x_{2},y_{2},...,x_{r},y_{r})).
\end{equation*}%
Let $\overrightarrow{P_{n}(z)}$ be a row vector of length $2^{2r}$ whose set
of components is%
\begin{equation*}
\{ P_{n}(z,x_{1},y_{1},x_{2},y_{2},...,x_{r},y_{r}):\forall i\in
\{1,2,...,r\},
\end{equation*}
\begin{equation*}
x_{i}=x_{i-1}-d_{i}^{x}\text{ and }%
y_{i}=y_{i-1}-d_{i}^{y}\text{ for some }d_{i}^{x},d_{i}^{y}\in
\{0,1\} \} .
\end{equation*}%
The order of these components in the vector is not important; we will need
to pick an order later to do matrix multiplication, but for now we will
leave this unspecified. The values of $\mathbf{R}$ at adjacent lattice
points can only differ by $1$ or $0$ so the vector $\overrightarrow{P_{n}(z)}
$ contains all possible values for vertices in the same section as $(n,n)$.
Thus%
\begin{equation*}
P_{n}(z)=\dsum\limits_{i=1}^{2^{2r}}\overrightarrow{P_{n}(z)}_{i}=%
\overrightarrow{P_{n}(z)}\mathbf{1}
\end{equation*}%
where $\mathbf{1}$ is the column vector $\mathbf{(}1,...,1)^{\prime}$.

Now we look at the relationship between $\overrightarrow{P_{n}(z)}$ and $%
\overrightarrow{P_{n-1}(z)}$. Let $x_{0}^{\prime }=y_{0}^{\prime }=z^{\prime
}$. If $\overrightarrow{P_{n}(z)}%
_{j}=P_{n}(z,x_{1},y_{1},x_{2},y_{2},...,x_{r},y_{r})$ and $\overrightarrow{%
P_{n-1}(z^{\prime })}_{i}=P_{n-1}(z^{\prime },x_{1}^{\prime },y_{1}^{\prime
},x_{2}^{\prime },y_{2}^{\prime },...,x_{r}^{\prime },y_{r}^{\prime })$ and $%
z^{\prime }\in \{z,$ $z-1\}$ define%
\begin{equation}
\Pr (R_{n}(z,x_{1},y_{1},x_{2},y_{2},...,x_{r},y_{r})\text{ and }%
R_{n-1}(z^{\prime },x_{1}^{\prime },y_{1}^{\prime },x_{2}^{\prime
},y_{2}^{\prime },...,x_{r}^{\prime },y_{r}^{\prime }))=\left\{
\begin{array}{ll}
\mathbf{M}_{ij} & \text{if }z^{\prime }=z \\
\mathbf{N}_{ij} & \text{if }z^{\prime }=z-1%
\end{array}%
\right.  \label{MN}
\end{equation}%
It sufficed to define this only for $z^{\prime }=z$ or $z-1$ because
otherwise the probability is $0$. Therefore summing over all possibilites
for $R_{n-1}()$ in the above expression gives us $\overrightarrow{P_{n}(z)}%
_{j}$:

\begin{equation*}
\dsum\limits_{i=1}\overrightarrow{P_{n-1}(z)}_{i}\mathbf{M}%
_{ij}+\dsum\limits_{i=1}\overrightarrow{P_{n-1}(z-1)}_{i}\mathbf{N}_{ij}=\Pr
(R_{n}(z,x_{1},y_{1},x_{2},y_{2},...,x_{r},y_{r}))=\overrightarrow{P_{n}(z)}%
_{j}
\end{equation*}%
Taking the convention that $\overrightarrow{P_{n}(z)}$ is the zero vector
for $n<r$, this yields the recurrence that is true for all $n\neq r$:

\begin{equation}
\overrightarrow{P_{n}(z)}=\overrightarrow{P_{n-1}(z)}\mathbf{M}+%
\overrightarrow{P_{n-1}(z-1)}\mathbf{N}\text{ }(n\neq r)  \label{recur}
\end{equation}

Now we will construct some generating functions. The convention made above
allows the generating function variables $n$ and $z$ to extend over all
integers. We will work with the two different generating functions $%
\overrightarrow{H_{n}(b)}=\dsum\limits_{z}$ $\overrightarrow{P_{n}(z)}b^{z}$
and $\overrightarrow{G(a,b)}=\dsum\limits_{n,z}$ $\overrightarrow{P_{n}(z)}%
a^{n}b^{z}$.

\subsection{The generating function $\protect\overrightarrow{G(a,b)}$}

Multiplying (\ref{recur}) by $a^{n}b^{z}$ and summing over all $n\neq r$ and
all $z$ yields%
\begin{equation*}
\dsum\limits_{n\neq r,z}\overrightarrow{P_{n}(z)}a^{n}b^{z}=\dsum\limits_{n%
\neq r,z}(\overrightarrow{P_{n-1}(z)}\mathbf{M}a^{n}b^{z})+\dsum\limits_{n%
\neq r,z}(\overrightarrow{P_{n-1}(z-1)}\mathbf{N}a^{n}b^{z})
\end{equation*}

Add $a^{r}\overrightarrow{H_{r}(b)}$ to both sides to obtain%
\begin{equation*}
\dsum\limits_{n,z}\overrightarrow{P_{n}(z)}a^{n}b^{z}=\left(
\dsum\limits_{n\neq r,z}\overrightarrow{P_{n-1}(z)}a^{n}b^{z}\right) \mathbf{%
M}+\left( \dsum\limits_{n\neq r,z}\overrightarrow{P_{n-1}(z-1)}%
a^{n}b^{z}\right) \mathbf{N}+a^{r}\overrightarrow{H_{r}(b)}
\end{equation*}%
Since $\overrightarrow{P_{r-1}(z)}$ is the zero vector, this becomes%
\begin{equation*}
\overrightarrow{G(a,b)}=a\overrightarrow{G(a,b)}\mathbf{M}+ab\overrightarrow{%
G(a,b)}\mathbf{N}+a^{r}\overrightarrow{H_{r}(b)}.
\end{equation*}

Then%
\begin{equation}
\overrightarrow{G(a,b)}(\mathbf{I}-a\mathbf{M}-ab\mathbf{N})=a^{r}%
\overrightarrow{H_{r}(b)}.  \label{twovarsgenfunction}
\end{equation}

\subsection{The generating function \protect\bigskip $\protect%
\overrightarrow{H_{n}(b)}$}

.We can also multiply (\ref{recur}) by $b^{z}$ and sum over all $z$ to obtain%
\begin{eqnarray*}
\dsum\limits_{z}\overrightarrow{P_{n}(z)}b^{z} &=&\left( \dsum\limits_{z}%
\overrightarrow{P_{n-1}(z)}b^{z}\right) \mathbf{M}+\left( \dsum\limits_{z}%
\overrightarrow{P_{n-1}(z-1)}b^{z}\right) \mathbf{N}\text{ \ }(n\neq r)\text{
}\Longrightarrow \\
\overrightarrow{H_{n}(b)} &=&\overrightarrow{H_{n-1}(b)}\mathbf{M+}b%
\overrightarrow{H_{n-1}(b)}\mathbf{N}\text{ \ }(n\neq r)
\end{eqnarray*}

This shows we can obtain $\overrightarrow{H_{n}(b)}$ by successive
multiplications by $\mathbf{M}+b\mathbf{N}$; that is, let $\mathbf{T(}b%
\mathbf{)}=\mathbf{M}+b\mathbf{N}$.\textbf{\ }%
\begin{equation*}
\overrightarrow{H_{n}(b)}=\overrightarrow{H_{r}(b)}\mathbf{T}(b)^{n-r}
\end{equation*}

To obtain the behavior of $\mathbf{T(}b\mathbf{)}^{n-r}$ as $%
n\longrightarrow \infty $ we assume from now on $b\geqslant 0$. We can then
apply results about positive matrices (see e.g. \cite{S}). Let $\det (%
\mathbf{T}(b)-\lambda \mathbf{I)=}g(\lambda ,b)$, a polynomial in $\lambda $
and $b$. $g(\lambda ,b)=(\lambda -f_{1}(b))(\lambda -f_{2}(b))...(\lambda
-f_{2^{2r}}(b))$. Let $\mathbf{e}(b)=(e_{1}(b),...,e_{r}(b))^{\prime }>%
\mathbf{0}$ be s.t. $\mathbf{T(}b\mathbf{)e}(b)=\mathbf{e}(b)f_{1}(b)$ and
let $\mathbf{e}^{\ast }(b)=(e_{1}^{\ast }(b),...,e_{r}^{\ast }(b))>\mathbf{0}%
^{\prime }$ s.t. $\mathbf{e}^{\ast }(b)f_{1}(b)=\mathbf{e}^{\ast }(b)\mathbf{%
T(}b\mathbf{)}$. Normalize\textbf{\ }$\mathbf{e}(b)$ and $\mathbf{e}^{%
\mathbf{\ast }}(b)$ so that $\mathbf{e}(b)\mathbf{1=}1,\mathbf{e}^{\ast }(b)%
\mathbf{1=}1$. Applying results for positive matrices,
\begin{equation}
\lim_{n\longrightarrow \infty }\frac{\mathbf{T(}b\mathbf{)}^{n}}{f_{1}(b%
\mathbf{)}^{n}}=\mathbf{e(}b\mathbf{)e}^{\mathbf{\ast }}(b)\Rightarrow
\lim_{n\longrightarrow \infty }\frac{(\mathbf{T(}b\mathbf{)}^{n})_{ij}}{%
nf_{1}(b\mathbf{)}^{n}}=0  \label{needed}
\end{equation}%
When $b=1$, this becomes%
\begin{equation}
\lim_{n\longrightarrow \infty }\mathbf{T(}1\mathbf{)}^{n}=\mathbf{1e}^{%
\mathbf{\ast }}(1)  \label{T(1)}
\end{equation}

since $\mathbf{T}(1)$ is the transition matrix between probability
distributions $\overrightarrow{H_{n-1}(1)}$ and $\overrightarrow{H_{n}(1)}.$

Let $h_{n}(b)=\frac{(\mathbf{T(}b\mathbf{)}^{n})_{ij}}{nf_{1}(b\mathbf{)}^{n}%
}$. We need the following limit result to complete the analysis. It appears
that it should follow from (\ref{needed}), but a proof eludes us. For now,
we will assume it to complete the analysis.

\begin{claim}
\begin{equation*}
\lim_{n\longrightarrow \infty }\frac{dh_{n}(b)}{db}=0
\end{equation*}
\end{claim}

The next step is
\begin{equation*}
\frac{dh_{n}(b)}{db}=\frac{1}{nf_{1}(b\mathbf{)}^{n}}\frac{d(\mathbf{T(}b%
\mathbf{)}^{n})_{ij}}{db}-\frac{(\mathbf{T(}b\mathbf{)}^{n})_{ij}}{f_{1}(b%
\mathbf{)}^{n+1}}\frac{df_{1}(b\mathbf{)}}{db}\text{ }\Rightarrow
\end{equation*}%
\begin{equation*}
\QOVERD. \vert {dh_{n}(b)}{db}_{b=1}=\frac{1}{n}\QOVERD. \vert {d(\mathbf{T(}%
b\mathbf{)}^{n})_{ij}}{db}_{b=1}-(\mathbf{T(}1\mathbf{)}^{n})_{ij}\QOVERD.
\vert {df_{1}(b\mathbf{)}}{db}_{b=1}\text{ }\Rightarrow
\end{equation*}%
\begin{equation}
\lim_{n\longrightarrow \infty }\frac{1}{n}\QOVERD. \vert {d(\mathbf{T(}b%
\mathbf{)}^{n})}{db}_{b=1}=\lim_{n\longrightarrow \infty }(\mathbf{T(}1%
\mathbf{)}^{n})\QOVERD. \vert {df_{1}(b\mathbf{)}}{db}_{b=1}\text{ =}\mathbf{%
1e}^{\mathbf{\ast }}(1)\QOVERD. \vert {df_{1}(b\mathbf{)}}{db}_{b=1}
\label{limresult}
\end{equation}%
Where the last implication follows from the unproven claim and (\ref{T(1)}).
Now we can apply this result to find $\mathbf{E}\overrightarrow{P_{n}(z)}$
which is defined below%
\begin{equation*}
\mathbf{E}\overrightarrow{P_{n}(z)}\equiv \dsum\limits_{z}z\overrightarrow{%
P_{n}(z)}=\QOVERD. \vert {d\overrightarrow{H_{n}(b)}}{db}_{b=1}
\end{equation*}%
Dividing by $n$ and taking the limit of both sides yields%
\begin{equation*}
\lim_{n\longrightarrow \infty }\frac{\mathbf{E}\overrightarrow{P_{n}(z)}}{n}%
=\lim_{n\longrightarrow \infty }\frac{1}{n}\QOVERD. \vert {d\overrightarrow{%
H_{n}(b)}}{db}_{b=1}=\lim_{n\longrightarrow \infty }\frac{1}{n}\QOVERD.
\vert {d(\overrightarrow{H_{r}(b)}\mathbf{T}(b)^{n-r}\mathbf{)}}{db}_{b=1}=
\end{equation*}%
\begin{equation*}
\lim_{n\longrightarrow \infty }\left( \frac{1}{n}\mathbf{T}(1)^{n-r}\QOVERD.
\vert {d(\overrightarrow{H_{r}(b)}\mathbf{)}}{db}_{b=1}+\frac{1}{n}%
\overrightarrow{H_{r}(1)}\QOVERD. \vert {d(\mathbf{T}(b\mathbf{)}^{n-r}%
\mathbf{)}}{db}_{b=1}\right) =
\end{equation*}%
\begin{equation*}
\lim_{n\longrightarrow \infty }\left( \frac{1}{n}\mathbf{T}(1)^{n-r}\mathbf{E%
}\overrightarrow{P_{r}(z)}+\frac{1}{n}\overrightarrow{H_{r}(1)}\QOVERD.
\vert {d(\mathbf{T}(b\mathbf{)}^{n-r}\mathbf{)}}{db}_{b=1}\right) =
\end{equation*}%
\begin{equation*}
\overrightarrow{H_{r}(1)}\lim_{n\longrightarrow \infty }\left( \frac{1}{n}%
\QOVERD. \vert {d(\mathbf{T}(b\mathbf{)}^{n-r}\mathbf{)}}{db}_{b=1}\right) =%
\overrightarrow{H_{r}(1)}\mathbf{1e}^{\mathbf{\ast }}(1)\QOVERD. \vert
{df_{1}(b\mathbf{)}}{db}_{b=1}=\mathbf{e}^{\mathbf{\ast }}(1)\QOVERD. \vert
{df_{1}(b\mathbf{)}}{db}_{b=1}
\end{equation*}

This last line uses (\ref{limresult}) and $\overrightarrow{H_{n}(1)}\mathbf{1%
}=\dsum\limits_{z}$ $\overrightarrow{P_{n}(z)}\mathbf{1}=\dsum%
\limits_{z}P_{n}(z)=1$. The equality above and the equation obtained by
multiplying it by $\mathbf{1}$ are stated below; they will be useful later.
\begin{equation}
\lim_{n\longrightarrow \infty }\frac{\mathbf{E}\overrightarrow{P_{n}(z)}}{n}=%
\mathbf{e}^{\mathbf{\ast }}(1)\QOVERD. \vert {df_{1}(b\mathbf{)}}{db}_{b=1}%
\text{ and }\lim_{n\longrightarrow \infty }\frac{\mathbf{EL}_{n,k,r}^{B}}{n}%
=\gamma _{k,r}^{B}=\QOVERD. \vert {df_{1}(b\mathbf{)}}{db}_{b=1}
\label{Epn(r)}
\end{equation}

The following claim makes computing $\QOVERD. \vert {df_{1}(b\mathbf{)}%
}{db}_{b=1}$ easier.

\begin{claim}
\label{rootclaim}Let $f_{1}(b)$ be the root of $g$ with $f_{1}(1)=1$. Then $%
\QOVERD. \vert {df_{1}(b)}{db}_{b=1}=-\QOVERD. \vert
{dg(1,b)}{db}_{b=1}\left. \QOVERD( ) {\lambda -1}{g(\lambda ,1)}\right\vert
_{\lambda =1}$
\end{claim}

\begin{proof}
\begin{equation*}
\frac{dg(1,b)}{db}=\frac{d(1-f_{1}(b))}{db}%
(1-f_{2}(b))...(1-f_{2^{2r}}(b))+(1-f_{1}(b))\frac{%
d((1-f_{2}(b))...(1-f_{2^{2r}}(b)))}{db}
\end{equation*}%
evaluating at $b=1$ yields%
\begin{equation*}
-\QOVERD. \vert {df_{1}(b\mathbf{)}%
}{db}_{b=1}(1-f_{2}(1))...(1-f_{2^{2r}}(1)).
\end{equation*}%
$g(\lambda ,1)=(\lambda -1)(\lambda -f_{2}(1))...(\lambda -f_{2^{2r}}(1))$
so $\lambda -1$ divides $g(\lambda ,1)$. $g(\lambda ,1)$ has only one root
at $\lambda =1$ because this root corresponds to the eigenvector $\mathbf{1}$
of $\mathbf{T}(1);\mathbf{1}$ is the unique positive eigenvector of $\mathbf{%
T}(1)$ (see e.g. \cite{S}). Thus $\frac{\lambda -1}{g(\lambda ,1)}$ is
defined at $\lambda =1$ and equals
\begin{equation*}
\frac{1}{(1-f_{2}(1))...(1-f_{2^{2r}}(1))}
\end{equation*}%
from which the claim follows directly.
\end{proof}

\subsection{Detailed analysis of 1-reach}

When $r=1$, $\overrightarrow{P_{n}(z)}%
=(P_{n}(z,z,z),P_{n}(z,z,z-1),P_{n}(z,z-1,z),P_{n}(z,z-1,z-1))$. The
matrices $\mathbf{M}$ and $\mathbf{N}$ are not difficult to compute by hand;
they are

\begin{eqnarray*}
\text{ \ \ \ \ \ }%
\begin{array}{l}
{\small P}_{n-1}{\small (z,z,z)} \\
{\small P}_{n-1}{\small (z,z,z-1)} \\
{\small P}_{n-1}{\small (z,z-1,z)} \\
{\small P}_{n-1}{\small (z,z-1,z-1)}%
\end{array}%
\begin{bmatrix}
\frac{(k-1)^{3}}{k^{3}} & 0 & 0 & 0 \\
\frac{(k-1)^{2}}{k^{2}} & 0 & 0 & 0 \\
\frac{(k-1)^{2}}{k^{2}} & 0 & 0 & 0 \\
\frac{k-1}{k} & 0 & 0 & 0%
\end{bmatrix}
&=&\mathbf{M}\text{ } \\
\text{\ }%
\begin{array}{c}
{\small P}_{n-1}{\small (z-1,z-1,z-1)} \\
{\small P}_{n-1}{\small (z-1,z-1,z-2)} \\
{\small P}_{n-1}{\small (z-1,z-2,z-1)} \\
{\small P}_{n-1}{\small (z-1,z-2,z-2)}%
\end{array}%
\begin{bmatrix}
\frac{1}{k^{2}} & \frac{k-1}{k^{2}} & \frac{k-1}{k^{2}} & \frac{(k-1)^{2}}{%
k^{3}} \\
0 & \frac{1}{k} & 0 & \frac{k-1}{k^{2}} \\
0 & 0 & \frac{1}{k} & \frac{k-1}{k^{2}} \\
0 & 0 & 0 & \frac{1}{k}%
\end{bmatrix}
&=&\mathbf{N}
\end{eqnarray*}

The expressions to the left of each matrix label the rows according to the
component order defined above; the columns correspond to $%
P_{n}(z,z,z),P_{n}(z,z,z-1),P_{n}(z,z-1,z),P_{n}(z,z-1,z-1)$ in that order.
We can also easily compute by hand $\overrightarrow{H_{1}(b)}=(\frac{k-1}{k}%
,0,0,\frac{b}{k})$.

\subsubsection{The two variable generating function}

(\ref{twovarsgenfunction}) gives us
\begin{equation*}
\overrightarrow{G(a,b)}=a\left( \frac{k-1}{k},0,0,\frac{b}{k}\right) (%
\mathbf{I}-a\mathbf{M}-ab\mathbf{N})^{-1}\text{.}
\end{equation*}%
Solving this problem with the two variable generating function is
computationally intensive, but it's nothing Maple can't handle. We obtain
\begin{equation*}
\overrightarrow{G(a,b)}^{\prime }=%
\begin{bmatrix}
ak^{2}(k-1)(-k+ab) \\
-a^{2}b(k-1)^{2}k \\
-a^{2}b(k-1)^{2}k \\
-ab(a^{2}b^{2}-abk^{2}-abk+k^{3})%
\end{bmatrix}%
\div
\end{equation*}%
\begin{equation*}
(a^{3}b^{3}-a^{2}b^{2}(k^{2}+2k)+a^{2}b(k^{3}-3k^{2}+3k-1)+ab(2k^{3}+k^{2})+a(k^{4}-3k^{3}+3k^{2}-k)-k^{4})
\end{equation*}%
This potentially gives us the entire distribution of the two variables. The
generating function for the expected value of $\overrightarrow{P_{n}(z)},$ $%
\mathbf{E}\overrightarrow{P_{n}(z)}=\dsum\limits_{z}z\overrightarrow{P_{n}(z)%
}$, is found by differentiating with respect to $b$ and then evaluating at $%
b=1$. We restrict to the $k=2$ case to make the expression simpler and more
readable.
\begin{equation*}
\dsum\limits_{n}\mathbf{E}\overrightarrow{P_{n}(z)}^{\prime }a^{n}=%
\begin{bmatrix}
-8a^{2}(a^{3}-7a^{2}+14a-12) \\
4a^{2}(a^{3}-4a^{2}-a+8) \\
4a^{2}(a^{3}-4a^{2}-a+8) \\
8a(3a^{3}-16a^{2}+26a-16)%
\end{bmatrix}%
(a^{3}-7a^{2}+22a-16)^{-2}
\end{equation*}

Using Mathematica's Discrete Math Rsolve package and a little computation by
hand, we get
\begin{equation*}
\mathbf{E}\overrightarrow{P_{n}(z)}^{\prime }=%
\begin{bmatrix}
\frac{32}{121}n-\frac{344}{1331}+2^{-2n}O(n) \\
\frac{16}{121}n-\frac{40}{1331}+2^{-2n}O(n) \\
\frac{16}{121}n-\frac{40}{1331}+2^{-2n}O(n) \\
\frac{24}{121}n+\frac{72}{1331}+2^{-2n}O(n)%
\end{bmatrix}%
\end{equation*}%
where the $O(n)$ terms vary like $n\cos (n\theta )$. Summing these
components gives us%
\begin{equation*}
\mathbf{EL}_{n,2,1}^{B}=\frac{8}{11}n-\frac{32}{121}+2^{-2n}O(n).
\end{equation*}

Mathematica can also solve the case for general $k$, but the expression is
difficult to pick apart because it's so long. To get the behavior of $\frac{%
\mathbf{EL}_{n,k,1}^{B}}{n}$ divide $\dsum\limits_{n}\mathbf{EL}%
_{n,k,1}^{B}a^{n}$ by $a$ and integrate with respect to $a$. This generating
function has the form
\begin{equation*}
\dsum\limits_{n}\frac{\mathbf{EL}_{n,k,1}^{B}}{n}a^{n}=\frac{c_{1}(k)}{1-a}%
+c_{2}(k)\ln (a-1)-c_{3}(k)\ln (O(a^{2}))+c_{4}(k)arctanh(O(a))
\end{equation*}

where $c_{i}(k)$ are functions only of $k$; the $O(a^{2})$ and $O(a)$ are
quadratic and linear polynomials in $a$ with coefficients a function of $k$.
Inferring from the $k=2$ case, we guess that%
\begin{equation*}
\mathbf{EL}_{n,k,1}^{B}=c_{1}(k)n-c_{2}(k)+2^{-2n}O(n)c_{5}(k).
\end{equation*}%
And Maple tells us that
\begin{equation*}
c_{1}(k)=\frac{3k+2}{(k^{2}+3k+1)}\text{, }c_{2}(k)=\frac{k(2k^{2}+3k+2)}{%
(k^{4}+6k^{3}+11k^{2}+6k+1)}
\end{equation*}

\subsubsection{The one variable generating function}

\begin{equation*}
\det (\mathbf{T}(b)-\lambda \mathbf{I)}\mathbf{=}g(\lambda ,b)=\frac{1}{k^{5}%
}(-\lambda k+b)\times
\end{equation*}%
\begin{equation}
(b^{3}-b^{2}k\lambda (k+2)+b\lambda (k^{3}+2k^{3}\lambda
-3k^{2}+k^{2}\lambda +3k-1)+\lambda ^{2}k(k^{3}-\lambda k^{3}-3k^{2}+3k-1)
\label{1reachg}
\end{equation}%
By (\ref{rootclaim})%
\begin{equation*}
\QOVERD. \vert {df_{1}(b)}{db}_{b=1}=-\QOVERD. \vert
{dg(1,b)}{db}_{b=1}\left. \QOVERD( ) {\lambda -1}{g(\lambda ,1)}\right\vert
_{\lambda =1}=
\end{equation*}%
\begin{equation*}
-\left( -\frac{1}{k^{5}}(k-1)^{3}(3k+2)\right) \left( \frac{k^{5}}{%
(k^{2}+3k+1)(k-1)^{3}}\right) =\frac{3k+2}{(k^{2}+3k+1)}.
\end{equation*}

Next we compute $\mathbf{e}^{\ast }(1)$ (using Maple even though it's not
necessary)

\begin{equation*}
\mathbf{e}^{\ast }(1)=N%
\begin{bmatrix}
k & 1 & 1 & \frac{1+k}{k}%
\end{bmatrix}%
\end{equation*}%
Choose $N$ so that $\mathbf{e}^{\ast }(1)\mathbf{1=}1$ $\Rightarrow N=\frac{k%
}{k^{2}+3k+1}$. From (\ref{Epn(r)}) we have

\begin{equation*}
\lim_{n\longrightarrow \infty }\frac{\mathbf{E}\overrightarrow{P_{n}(z)}}{n}=%
\mathbf{e}^{\mathbf{\ast }}(1)\QOVERD. \vert {df_{1}(b\mathbf{)}}{db}_{b=1}=n%
\frac{k(3k+2)}{(k^{2}+3k+1)^{2}}%
\begin{bmatrix}
k & 1 & 1 & \frac{1+k}{k}%
\end{bmatrix}%
\end{equation*}%
Summing all the components gives us%
\begin{equation*}
\gamma _{k,1}^{B}=\frac{3k+2}{(k^{2}+3k+1)}
\end{equation*}%
This does not give us as much asymptotic information as the two variable
generating function, but it is much less messy and allows us to easily see
the limiting behavior of $\mathbf{E}\overrightarrow{P_{n}(z)}$.

It is interesting to compare this limiting behavior to the conjectured
behavior for $\gamma _{k}^{B}$. It is guessed that $\sqrt{k}\gamma
_{k}^{B}\longrightarrow 2$ as $k\longrightarrow \infty $, whereas $k\gamma
_{k,1}^{B}\longrightarrow 3$ as $k\longrightarrow \infty $.

\subsection{2 and 3 reach}

When $r=2$, $\mathbf{M}$ and $\mathbf{N}$ are matrices of size $16\times 16$%
. For the two variable generating function approach, we will restrict to the
case $k=2$. Maple can solve for $\overrightarrow{G(a,b)}$; $\overrightarrow{%
G(a,b)}\mathbf{1}$ is a polynomial in $a$ and $b$ with leading term $%
a^{11}b^{11}$ divided by a polynomial with leading term $a^{10}b^{10}$. As
with 1-reach, we can find $\dint \left( \dsum\limits_{n}\mathbf{EL}%
_{n,2,2}^{B}a^{n-1}\right) da$ to obtain the limiting behavior of $\mathbf{EL%
}_{n,2,2}^{B}$. The result is an expression about a page long that is very
difficult to read. But it appears that most relevant parts of it to the
asymptotic behavior are:
\begin{equation*}
\frac{a(1-a)}{2(1-a)}+\frac{152}{197(1-a)}+\frac{16872(1-a)}{38809(1-a)}\ln
(a-1)
\end{equation*}%
From which we conclude
\begin{equation*}
\mathbf{EL}_{n,2,2}^{B}\mathbf{\sim }\frac{152}{197}n-\frac{16872}{38809}.
\end{equation*}%
This seems to be consistent with the Monte Carlo approximations, as will be
seen later.

Now for the one variable generating function approach. This can be solved
for general $k$. $g(\lambda ,b)$ is too large an expression to be of much
worth written down here. The resulting expression for $\QOVERD. \vert
{df_{1}(b)}{db}_{b=1}$ is surprisingly simple however.
\begin{equation*}
\QOVERD. \vert {df_{1}(b)}{db}_{b=1}=-\QOVERD. \vert
{dg(1,b)}{db}_{b=1}\left. \QOVERD( ) {\lambda -1}{g(\lambda ,1)}\right\vert
_{\lambda =1}=
\end{equation*}%
\begin{equation*}
-\left( -\frac{1}{k^{28}}%
(k+1)(5k^{3}+20k^{2}+15k+2)(k^{4}+k^{3}+3k^{2}+k+1)(k-1)^{15}(k^{4}+3k^{3}+5k^{2}+3k+1)\right) \times
\end{equation*}

\begin{equation*}
\left( \frac{k^{28}}{%
(k^{4}+3k^{3}+5k^{2}+3k+1)(k-1)^{15}(k+1)(k^{4}+k^{3}+3k^{2}+k+1)(k^{4}+10k^{3}+20k^{2}+10k+1)%
}\right) =
\end{equation*}%
\begin{equation*}
\frac{5k^{3}+20k^{2}+15k+2}{k^{4}+10k^{3}+20k^{2}+10k+1}=\gamma _{k,2}^{B}%
\text{.}
\end{equation*}%
when $k=2$, this gives $\frac{152}{197}$ which confirms part of the guess
for $\mathbf{EL}_{n,2,2}^{B}$ found by the two variable generating function
approach. $\mathbf{e}^{\mathbf{\ast }}(1)$ is illustrated as follows: We
reshape the vector into a matrix so that it is easier to read. The component
of $\mathbf{e}^{\mathbf{\ast }}(1)$ that corresponds to $%
P_{n}(z,z-d_{1}^{x},z-d_{1}^{y},z-d_{2}^{x},z-d_{2}^{y})$ in $%
\overrightarrow{P_{n}(z)}$ is represented by
\begin{tabular}{ll|l|}
\hline
\multicolumn{1}{|l}{$d_{2}^{x}$} & \multicolumn{1}{|l|}{$d_{1}^{x}$} & $0$
\\ \hline
&  & $d_{1}^{y}$ \\ \cline{3-3}
&  & $d_{2}^{y}$ \\ \cline{3-3}
\end{tabular}%
.
\begin{equation*}
\QATOP{%
\begin{bmatrix}
\begin{tabular}{ll|l|}
\hline
\multicolumn{1}{|l}{${\tiny 0}$} & \multicolumn{1}{|l|}{${\tiny 0}$} & $%
{\tiny 0}$ \\ \hline
&  & ${\tiny 0}$ \\ \cline{3-3}
&  & ${\tiny 0}$ \\ \cline{3-3}
\end{tabular}
&
\begin{tabular}{ll|l|}
\hline
\multicolumn{1}{|l}{${\tiny 0}$ \ } & \multicolumn{1}{|l|}{${\tiny 0}$ \ } &
${\tiny 0}$ \\ \hline
&  & ${\tiny 0}$ \\ \cline{3-3}
&  & ${\tiny 1}$ \\ \cline{3-3}
\end{tabular}
&
\begin{tabular}{ll|l|}
\hline
\multicolumn{1}{|l}{${\tiny 0}$ \ } & \multicolumn{1}{|l|}{${\tiny 0}$ \ } &
${\tiny 0}$ \\ \hline
&  & ${\tiny 1}$ \\ \cline{3-3}
&  & ${\tiny 1}$ \\ \cline{3-3}
\end{tabular}
&
\begin{tabular}{ll|l|}
\hline
\multicolumn{1}{|l}{${\tiny 0}$ \ } & \multicolumn{1}{|l|}{${\tiny 0}$ \ } &
${\tiny 0}$ \\ \hline
&  & ${\tiny 1}$ \\ \cline{3-3}
&  & ${\tiny 2}$ \\ \cline{3-3}
\end{tabular}
\\
\begin{tabular}{ll|l|}
\hline
\multicolumn{1}{|l}{${\tiny 1}$} & \multicolumn{1}{|l|}{${\tiny 0}$ \ } & $%
{\tiny 0}$ \\ \hline
&  & ${\tiny 0}$ \\ \cline{3-3}
&  & ${\tiny 0}$ \\ \cline{3-3}
\end{tabular}
&
\begin{tabular}{ll|l|}
\hline
\multicolumn{1}{|l}{${\tiny 1}$} & \multicolumn{1}{|l|}{${\tiny 0}$ \ } & $%
{\tiny 0}$ \\ \hline
&  & ${\tiny 0}$ \\ \cline{3-3}
&  & ${\tiny 1}$ \\ \cline{3-3}
\end{tabular}
&
\begin{tabular}{ll|l|}
\hline
\multicolumn{1}{|l}{${\tiny 1}$} & \multicolumn{1}{|l|}{${\tiny 0}$ \ } & $%
{\tiny 0}$ \\ \hline
&  & ${\tiny 1}$ \\ \cline{3-3}
&  & ${\tiny 1}$ \\ \cline{3-3}
\end{tabular}
&
\begin{tabular}{ll|l|}
\hline
\multicolumn{1}{|l}{${\tiny 1}$} & \multicolumn{1}{|l|}{${\tiny 0}$ \ } & $%
{\tiny 0}$ \\ \hline
&  & ${\tiny 1}$ \\ \cline{3-3}
&  & ${\tiny 2}$ \\ \cline{3-3}
\end{tabular}
\\
\begin{tabular}{ll|l|}
\hline
\multicolumn{1}{|l}{${\tiny 1}$} & \multicolumn{1}{|l|}{${\tiny 1}$} & $%
{\tiny 0}$ \\ \hline
&  & ${\tiny 0}$ \\ \cline{3-3}
&  & ${\tiny 0}$ \\ \cline{3-3}
\end{tabular}
&
\begin{tabular}{ll|l|}
\hline
\multicolumn{1}{|l}{${\tiny 1}$} & \multicolumn{1}{|l|}{${\tiny 1}$} & $%
{\tiny 0}$ \\ \hline
&  & ${\tiny 0}$ \\ \cline{3-3}
&  & ${\tiny 1}$ \\ \cline{3-3}
\end{tabular}
&
\begin{tabular}{ll|l|}
\hline
\multicolumn{1}{|l}{${\tiny 1}$} & \multicolumn{1}{|l|}{${\tiny 1}$} & $%
{\tiny 0}$ \\ \hline
&  & ${\tiny 1}$ \\ \cline{3-3}
&  & ${\tiny 1}$ \\ \cline{3-3}
\end{tabular}
&
\begin{tabular}{ll|l|}
\hline
\multicolumn{1}{|l}{${\tiny 1}$} & \multicolumn{1}{|l|}{${\tiny 1}$} & $%
{\tiny 0}$ \\ \hline
&  & ${\tiny 1}$ \\ \cline{3-3}
&  & ${\tiny 2}$ \\ \cline{3-3}
\end{tabular}
\\
\begin{tabular}{ll|l|}
\hline
\multicolumn{1}{|l}{${\tiny 2}$} & \multicolumn{1}{|l|}{${\tiny 1}$} & $%
{\tiny 0}$ \\ \hline
&  & ${\tiny 0}$ \\ \cline{3-3}
&  & ${\tiny 0}$ \\ \cline{3-3}
\end{tabular}
&
\begin{tabular}{ll|l|}
\hline
\multicolumn{1}{|l}{${\tiny 2}$} & \multicolumn{1}{|l|}{${\tiny 1}$} & $%
{\tiny 0}$ \\ \hline
&  & ${\tiny 0}$ \\ \cline{3-3}
&  & ${\tiny 1}$ \\ \cline{3-3}
\end{tabular}
&
\begin{tabular}{ll|l|}
\hline
\multicolumn{1}{|l}{${\tiny 2}$} & \multicolumn{1}{|l|}{${\tiny 1}$} & $%
{\tiny 0}$ \\ \hline
&  & ${\tiny 1}$ \\ \cline{3-3}
&  & ${\tiny 1}$ \\ \cline{3-3}
\end{tabular}
&
\begin{tabular}{ll|l|}
\hline
\multicolumn{1}{|l}{${\tiny 2}$} & \multicolumn{1}{|l|}{${\tiny 1}$} & $%
{\tiny 0}$ \\ \hline
&  & ${\tiny 1}$ \\ \cline{3-3}
&  & ${\tiny 2}$ \\ \cline{3-3}
\end{tabular}%
\end{bmatrix}%
}{\Updownarrow }
\end{equation*}%
\begin{equation*}
\begin{bmatrix}
{k}^{2} & 2\,{k} & {k} & 1 \\
2k & k+4 & k+2 & \frac{2\,\left( k+1\right) }{k} \\
{k} & k+2 & k+1 & \frac{2k+1}{k} \\
1 & \frac{2\,\left( k+1\right) }{k} & \frac{2k+1}{k} & \frac{{k}^{2}+4\,k+1}{%
k^{2}}%
\end{bmatrix}%
\end{equation*}%
To normalize $\mathbf{e}^{\mathbf{\ast }}(1)$, the above matrix must be
multiplied by $\frac{k^{2}}{k^{4}+10k^{3}+20k^{2}+10k+1}$ which finally gives%
\begin{equation*}
\mathbf{E}\overrightarrow{P_{n}(z)}\text{ }\mathbf{\sim }\text{ }n\frac{%
k^{2}(5k^{3}+20k^{2}+15k+2)}{(k^{4}+10k^{3}+20k^{2}+10k+1)^{2}}%
\begin{bmatrix}
{k}^{2} & 2\,{k} & {k} & 1 \\
2k & k+4 & k+2 & \frac{2\,\left( k+1\right) }{k} \\
{k} & k+2 & k+1 & \frac{2k+1}{k} \\
1 & \frac{2\,\left( k+1\right) }{k} & \frac{2k+1}{k} & \frac{{k}^{2}+4\,k+1}{%
k^{2}}%
\end{bmatrix}%
\end{equation*}

The case $r=3,$ $k=2$ is also computable in a reasonable amount of time (it
took Maple about a half an hour on a 1992 Mega Hertz Dell). The result is

\begin{equation*}
\gamma _{2,3}^{B}=\frac{3376}{4279}\text{.}
\end{equation*}

\section{Applications to the Random String model}

The machinery developed for r-reach with the Bernoulli matching model can be
applied to 1-reach with the Random String model when $k=2$. For $r>1$, it
appears this same brute force conditional probability approach is so
complicated as to be almost useless. $r=1$ and $k>2$ seems significantly
more difficult than $r=1$, $k=2$, which is rather surprising. We get an
interesting reduction for the $k=2$ case, as will be seen shortly. The
reason for pursuing this approach despite its appearance of being difficult
to generalize, is that it may lead to a short proof of $\gamma
_{2,1}^{B}>\gamma _{2,1}$, which may be generalizable. It has been
conjectured that $\lim_{n\longrightarrow \infty }\gamma _{k}^{B}\sqrt{k}%
=\lim_{n\longrightarrow \infty }\gamma _{k}\sqrt{k}$ (actually Sankoff and
Mainville conjectured that $\lim_{n\longrightarrow \infty }\gamma _{k}\sqrt{k%
}=2$ (see e.g. \cite{D}) and Boutet de Monvel \cite{B} conjectured that $%
\lim_{n\longrightarrow \infty }\gamma _{k}^{B}\sqrt{k}=2$). If 1-reach is
solved for general $k$, it may provide some insights into this problem.

\subsection{Detailed analysis of 1-reach}

The reduction for the case $k=2$ is not difficult, but it requires a fair
amount of notation to discuss.

\begin{definition}
If $\epsilon _{ij}$ is defined for $|i-j|\leq r$ and $1\leq i,j\leq n$, $%
\epsilon _{ij}$ is a string realizable configuration of weight $w$ if $%
\epsilon _{ij}=\delta _{u(i),v(j)}$ for $w$ distinct $(u,v)\epsilon \Sigma
^{n}x\Sigma ^{n}.$
\end{definition}

It is easy to convince oneself of the following claim by doing a case by
case analysis for $n=3$. Such an analysis extends easily to general $n$.

\begin{claim}
Let $k=2$ and let $\epsilon _{ij}$ be defined for $|i-j|\leq 1$ and $1\leq
i,j\leq n$. $\epsilon _{ij}$ is a string realizable configuration of weight $%
2$ if%
\begin{equation}
\forall i\in \{1,...,n\}\text{, }\epsilon _{i-1,i-1}+\epsilon
_{i,i-1}+\epsilon _{i-1,i\text{ }}+\epsilon _{i,i}\in \{0,2,4\}
\label{prop024}
\end{equation}%
and is a string realizable configuration of weight $0$ otherwise.
\end{claim}

\begin{proof}
$k=2$ means the alphabet, $\Sigma $, is $\{0,1\}$ so that $%
\begin{bmatrix}
\delta _{u(i-1),v(i)} & \delta _{u(i),v(i)} \\
\delta _{u(i-1),v(i-1)} & \delta _{u(i),v(i-1)}%
\end{bmatrix}%
\equiv \mathbf{X(}u(i-1,i),v(i-1,i))$ must be in the set
\begin{equation}
\mathbf{Y}\equiv \left\{
\begin{bmatrix}
1 & 1 \\
1 & 1%
\end{bmatrix}%
,%
\begin{bmatrix}
0 & 1 \\
1 & 0%
\end{bmatrix}%
,%
\begin{bmatrix}
0 & 0 \\
1 & 1%
\end{bmatrix}%
,%
\begin{bmatrix}
1 & 0 \\
1 & 0%
\end{bmatrix}%
,%
\begin{bmatrix}
0 & 0 \\
0 & 0%
\end{bmatrix}%
,%
\begin{bmatrix}
1 & 0 \\
0 & 1%
\end{bmatrix}%
,%
\begin{bmatrix}
1 & 1 \\
0 & 0%
\end{bmatrix}%
,%
\begin{bmatrix}
0 & 1 \\
0 & 1%
\end{bmatrix}%
\right\} .  \label{8config}
\end{equation}%
This shows that if the condition in (\ref{prop024}) fails, $\epsilon _{ij}$
is a string realizable configuration of weight $0$.

For the other part of the claim we proceed by induction on $n$. The case $n=1
$ can be seen by noting that each element of $\mathbf{Y}$ is equal to $2$ of
the $16$ possibilities for $\mathbf{X}(u(i-1,i),v(i-1,i))$. Suppose $n>1$
and the claim holds for $n-1$. Let $u,v,u^{\prime },v^{\prime }$ be the
strings of length $n-1$ such that $\forall i,j\in \{1,...,n-1\}$ and $%
|i-j|\leq 1,$ $\epsilon _{ij}=\delta _{u(i),v(j)}=\delta _{u^{\prime
}(i),v^{\prime }(j)}$. By hypothesis, $\mathbf{Z}_{n}\equiv
\begin{bmatrix}
\epsilon _{n-1,n} & \epsilon _{nn} \\
\epsilon _{n-1,n-1} & \epsilon _{n,n-1}%
\end{bmatrix}%
$ is one of the eight matrices belonging to $\mathbf{Y}$. For each matrix in
$\mathbf{Y}$, we can choose $u(n)$ and $v(n)$ as shown below so that $%
\forall i,j\in \{1,...,n\}$ and $|i-j|\leq 1,$ $\epsilon _{ij}=\delta
_{u(i),v(j)}$. The same goes for $u^{\prime }$ and $v^{\prime }$.
\begin{equation*}
\begin{tabular}{l|llllllll}
$\mathbf{Z}_{n}$ &
\begin{tabular}{|l|l|}
\hline
{\small 1} & {\small 1} \\ \hline
{\small 1} & {\small 1} \\ \hline
\end{tabular}
&
\begin{tabular}{|l|l|}
\hline
{\small 0} & {\small 1} \\ \hline
{\small 1} & {\small 0} \\ \hline
\end{tabular}
&
\begin{tabular}{|l|l|}
\hline
{\small 0} & {\small 0} \\ \hline
{\small 1} & {\small 1} \\ \hline
\end{tabular}
&
\begin{tabular}{|l|l|}
\hline
{\small 1} & {\small 0} \\ \hline
{\small 1} & {\small 0} \\ \hline
\end{tabular}
&
\begin{tabular}{|l|l|}
\hline
{\small 0} & {\small 0} \\ \hline
{\small 0} & {\small 0} \\ \hline
\end{tabular}
&
\begin{tabular}{|l|l|}
\hline
{\small 1} & {\small 0} \\ \hline
{\small 0} & {\small 1} \\ \hline
\end{tabular}
&
\begin{tabular}{|l|l|}
\hline
{\small 1} & {\small 1} \\ \hline
{\small 0} & {\small 0} \\ \hline
\end{tabular}
&
\begin{tabular}{|l|l|}
\hline
{\small 0} & {\small 1} \\ \hline
{\small 0} & {\small 1} \\ \hline
\end{tabular}
\\ \hline
$\QATOP{u(n)=,}{v(n)=}$ & $\QATOP{u(n-1),}{v(n-1)}$ & $\QATOP{1-u(n-1),}{%
1-v(n-1)}$ & $\QATOP{1-u(n-1),}{v(n-1)}$ & $\QATOP{u(n-1),}{1-v(n-1)}$ & $%
\QATOP{u(n-1),}{v(n-1)}$ & $\QATOP{1-u(n-1),}{1-v(n-1)}$ & $\QATOP{1-u(n-1),%
}{v(n-1)}$ & $\QATOP{u(n-1),}{1-v(n-1)}$%
\end{tabular}%
\end{equation*}%
This shows $\epsilon _{ij}$ is a string realizable configuration of weight
at least $2$. The weight cannot exceed $2$ because then $\epsilon _{ij}$
restricted to $i,j\in \{1,...,n-1\}$ would have weight greater than $2$.
\end{proof}

This claim lets us count the probabilities $P_{n}(z,x_{1},y_{1})$ much like
we did for the Bernoulli Matching model. We define the analogous probability
vector but we have to break $P_{n}(z,x_{1},y_{1})$ into two pieces: $%
P_{n}(z,x_{1},y_{1})=$ $P_{n}^{\text{on}}(z,x_{1},y_{1})+P_{n}^{\text{off}%
}(z,x_{1},y_{1})$.%
\begin{equation*}
P_{n}^{\text{on}}(z,x_{1},y_{1})=\Pr (R_{n}(z,x_{1},y_{1})\text{ and }%
\epsilon _{nn}=1)\text{, }P_{n}^{\text{off}}(z,x_{1},y_{1})=\Pr
(R_{n}(z,x_{1},y_{1})\text{ and }\epsilon _{nn}=0).
\end{equation*}%
The reason for this split is that we need to know $\epsilon _{n-1,n-1}$ to
determine how $\{\mathbf{R}_{n-1,n-1},\mathbf{R}_{n-2,n-1},\mathbf{R}%
_{n-1,n-2}\}$ affects $\{\mathbf{R}_{n,n},\mathbf{R}_{n-1,n},\mathbf{R}%
_{n,n-1}\}$. The computation of $\mathbf{M}$ and $\mathbf{N}$ was done by
hand and was a little trickier than for the Bernoulli Matching model.

\begin{equation*}
\text{ \ \ \ \ }%
\begin{array}{l}
{\small P}_{n-1}^{\text{off}}{\small (z,z,z)} \\
{\small P}_{n-1}^{\text{off}}{\small (z,z,z-1)} \\
{\small P}_{n-1}^{\text{off}}{\small (z,z-1,z)} \\
{\small P}_{n-1}^{\text{off}}{\small (z,z-1,z-1)} \\
{\small P}_{n-1}^{\text{on}}{\small (z,z,z)} \\
{\small P}_{n-1}^{\text{on}}{\small (z,z,z-1)} \\
{\small P}_{n-1}^{\text{on}}{\small (z,z-1,z)} \\
{\small P}_{n-1}^{\text{on}}{\small (z,z-1,z-1)}%
\end{array}%
\begin{bmatrix}
{\small 1/4} & {\small 0} & {\small 0} & {\small 0} & {\small 0} & {\small 0}
& {\small 0} & {\small 0} \\
{\small 1/4} & {\small 0} & {\small 0} & {\small 0} & {\small 0} & {\small 0}
& {\small 0} & {\small 0} \\
{\small 1/4} & {\small 0} & {\small 0} & {\small 0} & {\small 0} & {\small 0}
& {\small 0} & {\small 0} \\
{\small 1/2} & {\small 0} & {\small 0} & {\small 0} & {\small 0} & {\small 0}
& {\small 0} & {\small 0} \\
{\small 0} & {\small 0} & {\small 0} & {\small 0} & {\small 0} & {\small 0}
& {\small 0} & {\small 0} \\
{\small 1/4} & {\small 0} & {\small 0} & {\small 0} & {\small 0} & {\small 0}
& {\small 0} & {\small 0} \\
{\small 1/4} & {\small 0} & {\small 0} & {\small 0} & {\small 0} & {\small 0}
& {\small 0} & {\small 0} \\
{\small 1/2} & {\small 0} & {\small 0} & {\small 0} & {\small 0} & {\small 0}
& {\small 0} & {\small 0}%
\end{bmatrix}%
=\mathbf{M}
\end{equation*}%
\begin{equation*}
\text{ \ \ \ \ }%
\begin{array}{l}
{\small P}_{n-1}^{\text{off}}{\small (z-1,z-1,z-1)} \\
{\small P}_{n-1}^{\text{off}}{\small (z-1,z-1,z-2)} \\
{\small P}_{n-1}^{\text{off}}{\small (z-1,z-2,z-1)} \\
{\small P}_{n-1}^{\text{off}}{\small (z-1,z-2,z-2)} \\
{\small P}_{n-1}^{\text{on}}{\small (z-1,z-1,z-1)} \\
{\small P}_{n-1}^{\text{on}}{\small (z-1,z-1,z-2)} \\
{\small P}_{n-1}^{\text{on}}{\small (z-1,z-2,z-1)} \\
{\small P}_{n-1}^{\text{on}}{\small (z-1,z-2,z-2)}%
\end{array}%
\begin{bmatrix}
{\small 1/4} & {\small 0} & {\small 0} & {\small 0} & {\small 0} & {\small %
1/4} & {\small 1/4} & {\small 0} \\
{\small 0} & {\small 1/4} & {\small 0} & {\small 0} & {\small 0} & {\small %
1/4} & {\small 0} & {\small 1/4} \\
{\small 0} & {\small 0} & {\small 1/4} & {\small 0} & {\small 0} & {\small 0}
& 1{\small /4} & {\small 1/4} \\
{\small 0} & {\small 0} & {\small 0} & {\small 0} & {\small 0} & {\small 0}
& {\small 0} & {\small 1/2} \\
{\small 0} & 1/4 & {\small 1/4} & {\small 0} & {\small 1/4} & {\small 0} &
{\small 0} & {\small 1/4} \\
{\small 0} & 1/4 & {\small 0} & {\small 0} & {\small 0} & 1/4 & {\small 0} &
{\small 1/4} \\
{\small 0} & {\small 0} & {\small 1/4} & {\small 0} & {\small 0} & {\small 0}
& {\small 1/4} & {\small 1/4} \\
{\small 0} & {\small 0} & {\small 0} & {\small 0} & {\small 0} & {\small 0}
& {\small 0} & 1/2%
\end{bmatrix}%
=\mathbf{N}
\end{equation*}%
The two variable generating function approach determines the fine limiting
behavior:
\begin{equation*}
\dsum\limits_{n}\mathbf{EL}_{n,2,1}a^{n}=\frac{a(a^{2}-2a+8)}{%
2(a^{2}-4a+8)(a-1)^{2}}
\end{equation*}%
Using Mathematica's Discrete Math Rsolve package and a little computation by
hand, we obtain%
\begin{equation*}
\mathbf{EL}_{n,2,1}=\frac{7}{10}n-\frac{7}{25}+2^{-3n/2}O(1).
\end{equation*}%
where the $O(1)$ term varies like $\cos (n\theta )$. We will compare this
result to numerical approximations.

The one variable generating function produces

\begin{equation*}
\det (\mathbf{T}(b)-\lambda \mathbf{I)}\mathbf{=}g(\lambda ,b)=-\frac{1}{128}%
\lambda ^{3}(b-2\lambda )(b-4\lambda )(b^{3}+b^{2}(-8\lambda +1)+2b\lambda
(-1+10\lambda )+4\lambda ^{2}(-4\lambda +1))
\end{equation*}%
and%
\begin{equation*}
\QOVERD. \vert {df_{1}(b)}{db}_{b=1}=-\QOVERD. \vert
{dg(1,b)}{db}_{b=1}\left. \QOVERD( ) {\lambda -1}{g(\lambda ,1)}\right\vert
_{\lambda =1}=
\end{equation*}%
\begin{equation*}
-\left( -\frac{21}{128}\right) \left( \frac{64}{15}\right) =\frac{7}{10}
\end{equation*}%
$\mathbf{e}^{\mathbf{\ast }}(1)=\frac{1}{20}%
\begin{bmatrix}
8 & 1 & 1 & 0 & 0 & 3 & 3 & 4%
\end{bmatrix}%
$.

We can also "blow up" the 1-reach Bernoulli Matching model, so that we work
with $\overrightarrow{P_{n}^{\text{on}}(z)}$ and $\overrightarrow{P_{n}^{%
\text{off}}(z)}$ even though we don't need to. The resulting matrices are
included in the appendix. It is interesting to note that the matrices only
differ in the two rows corresponding to $P_{n-1}(z-1,z-1,z-1)$. The result is%
\begin{equation*}
g^{B}(\lambda ,b)=\frac{1}{32}\lambda ^{4}(b-2\lambda )(b^{3}-8b^{2}\lambda
+b\lambda (1+20\lambda )+2\lambda ^{2}(1-8\lambda ))
\end{equation*}%
and this polynomial is the same as one obtained earlier (in (\ref{1reachg}))
except for the $\lambda ^{4}$ term. Also,

$\mathbf{e}^{B\mathbf{\ast }}(1)=\frac{1}{22}%
\begin{bmatrix}
7 & 2 & 2 & 0 & 1 & 2 & 2 & 6%
\end{bmatrix}%
$ which is more precise behavior than that determined by the $4\times 4$
matrix method.

It is unclear whether there is a more direct way to see that the difference
in the matrices for the Random String model and the Bernoulli Matching model
lead to the conclusion $\gamma _{2,1}<$ $\gamma _{2,1}^{B}$

\section{Numerical Work}

We ran Monte Carlo simulations for $k=2$ and $r=1,2...10,15,20,25,35,40$. $%
10000$ trials were computed up to $n=1000$ for each $r$. To obtain behavior
varying with $n$, approximations of $\mathbf{EL}_{n,2,r}^{B}$ and $\mathbf{EL%
}_{n,2,r}$ for all $n$ from $1$ to $1000$ were computed for each trial.
Ideally, we should have computed separate trials for each $n$, but these
results appear to lead to good extrapolations to large n. Following the work
in \cite{B}, we extrapolate to $\gamma _{2,r}^{B}$ from the small n
simulations based on;
\begin{equation}
\mathbf{EL}_{n,2,r}\text{ }\mathbf{\sim }\text{ }\gamma _{2,r}n-A_{r},\text{
}\mathbf{EL}_{n,2,r}^{B}\text{ }\mathbf{\sim }\text{ }\gamma
_{2,r}^{B}n-A_{r}^{B}.  \label{finitesizefix}
\end{equation}

\begin{figure}[tbp]
\begin{center}
\includegraphics [scale=.65]{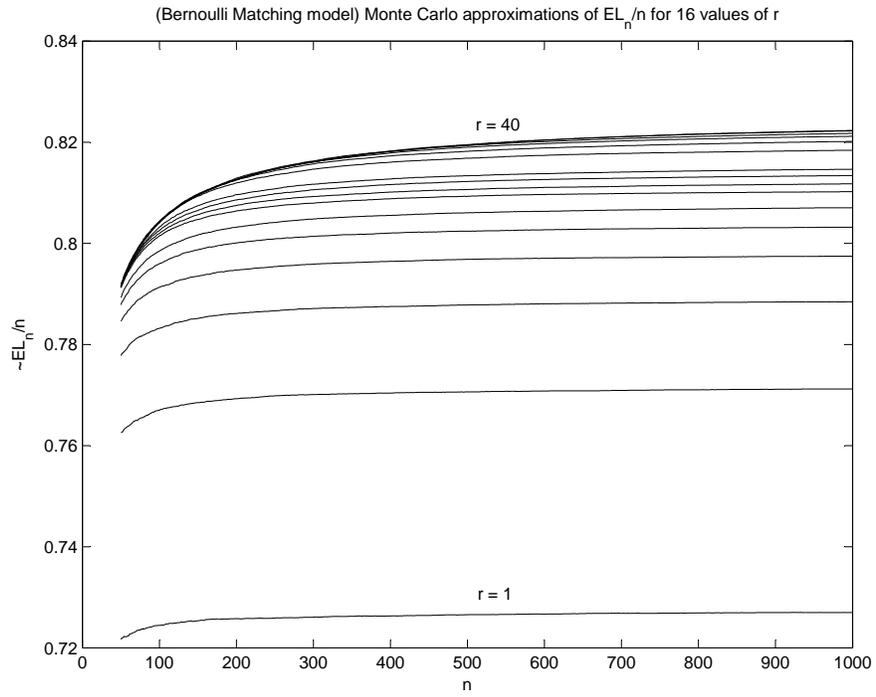}%
\vspace{0.4in}
\includegraphics[scale=.65]
{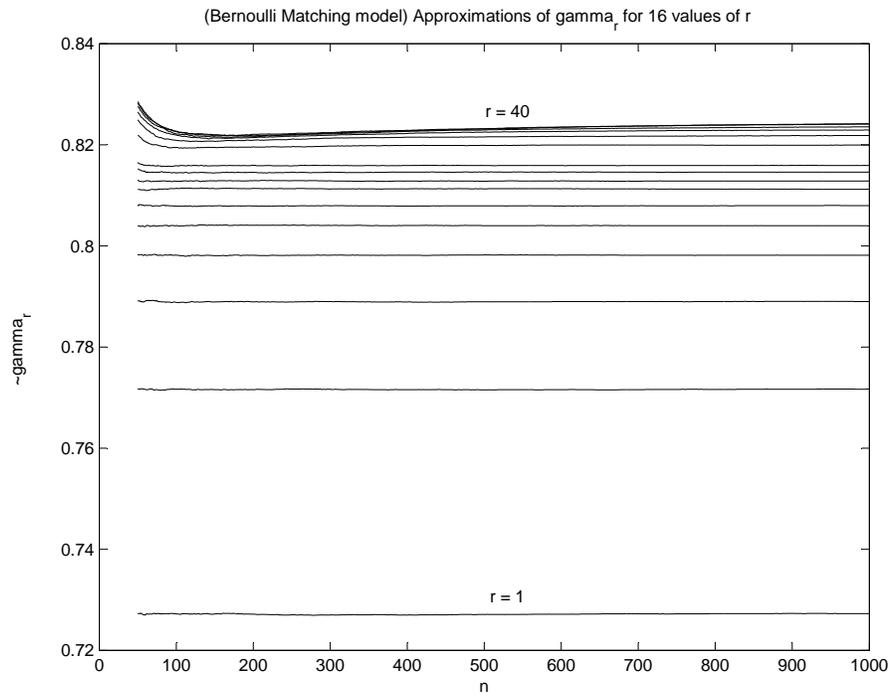}
\end{center}
\caption{The Monte Carlo approximations of $\frac{\mathbf{EL}_{n,2,r}^{B}}{n}
$ and this same data corrected by (\protect\ref{finitesizefix})\ to obtain
the limiting behavior.}
\label{twographsB}
\end{figure}

Where $A_{r}$ $(A_{r}^{B})$ is a constant, and was found by minimizing the
variance of $\frac{\mathbf{EL}_{n,2,r}\mathbf{+}A_{r}}{n}$ $(\frac{\mathbf{EL%
}_{n,2,r}^{B}\mathbf{+}A_{r}^{B}}{n})$. Extrapolations for $\gamma
_{2,r}^{B} $, $A_{r}^{B}$, and $\gamma _{2,r}$, $A_{r}$ based on Monte Carlo
simulations are shown below. We did this extrapolation from $n=50...1000$ to
minimize the effect of the $2^{-2n}O(n)$ term (we only saw this for $r=1$,
but there are probably similar terms for larger $r$).

\begin{equation*}
\begin{tabular}{l|lllllllll}
$r$ & $1$ & $2$ & $3$ & $4$ & $5$ & $6$ & $7$ & $8$ & $9$  \\ \hline
$\gamma _{2,r}^{B}$ & $0.72726$ & $0.77166$ & $0.78898$ & $0.79813$ & $%
0.80396$ & $0.80796$ & $0.81119$ & $0.81284$ & $0.81458$  \\
\hline
$A_{r}^{B}$ & $0.2771$ & $0.4626$ & $0.5641$ & $0.6852$ & $0.8033$ & $0.9399$
& $0.9931$ & $1.0814$ & $1.1900$  \\ \hline
$\gamma _{2,r}$ & $0.70014$ & $0.73767$ & $0.75610$ & $0.76718$ & $0.77467$
& $0.78004$ & $0.78408$ & $0.78726$ & $0.78976$  \\ \hline
$A_{r}$ & $0.2652$ & $0.4335$ & $0.5748$ & $0.7048$ & $0.8195$ & $0.9218$ & $%
1.0163$ & $1.1121$ & $1.2044$ %
\end{tabular}%
\end{equation*}%
\begin{equation*}
\begin{tabular}{l|lllllll}
$r$ & $10$ & $15$ & $20$ & $25$ & $30$ & $35$ & $40$ \\ \hline
$\gamma _{2,r}^{B}$ & $0.81592$ & $0.81994$ & $0.82182$ & $0.82290$ & $0.82355$ & $%
0.82406$ & $0.82415$ \\ \hline
$A_{r}^{B}$ & $1.2653$ & $1.5253$ & $1.6814$ & $1.7536$ & $1.8058$ & $1.8368$ & $1.8395$
\\ \hline
$\gamma _{2,r}$ & $0.79180$ & $0.79819$ & $0.80149$ & $0.80340$ & $0.80462$ & $0.80546$
& $0.80603$ \\ \hline
$A_{r}$ & $1.2877$ & $1.6377$ & $1.8753$ & $2.028$ & $2.1273$ & $2.1939$ & $2.2371$%
\end{tabular}%
\end{equation*}%
Shown in figure (\ref{twographsB}) are $\frac{MonteCarlo(\mathbf{EL}%
_{n,2,r}^{B})}{n}$ and $\frac{MonteCarlo(\mathbf{EL}_{n,2,r}^{B})\mathbf{+}%
A_{r}}{n}$ and the corresponding Random String model data is shown in (\ref%
{twographs}). It appears that the approximation $\mathbf{EL}_{n,2,r}$ $%
\mathbf{\sim }$ $\gamma _{2,r}n-A_{r}$ gets increasingly worse for larger $r$
and likewise for the Bernoulli Matching model.

\begin{figure}[tbp]
\begin{center}
\includegraphics [scale=.65]{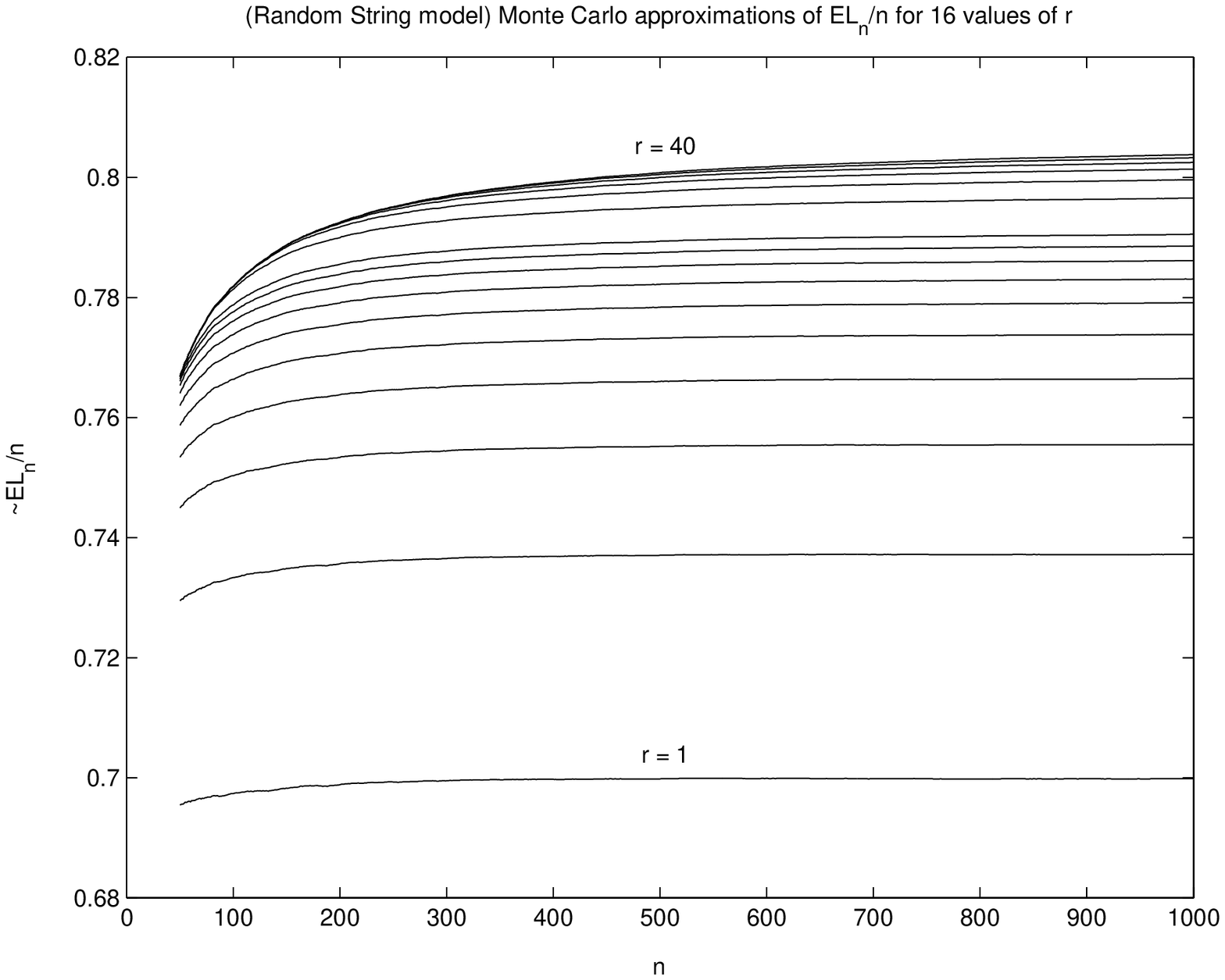}%
\vspace{0.4in}
\includegraphics[scale=.65]
{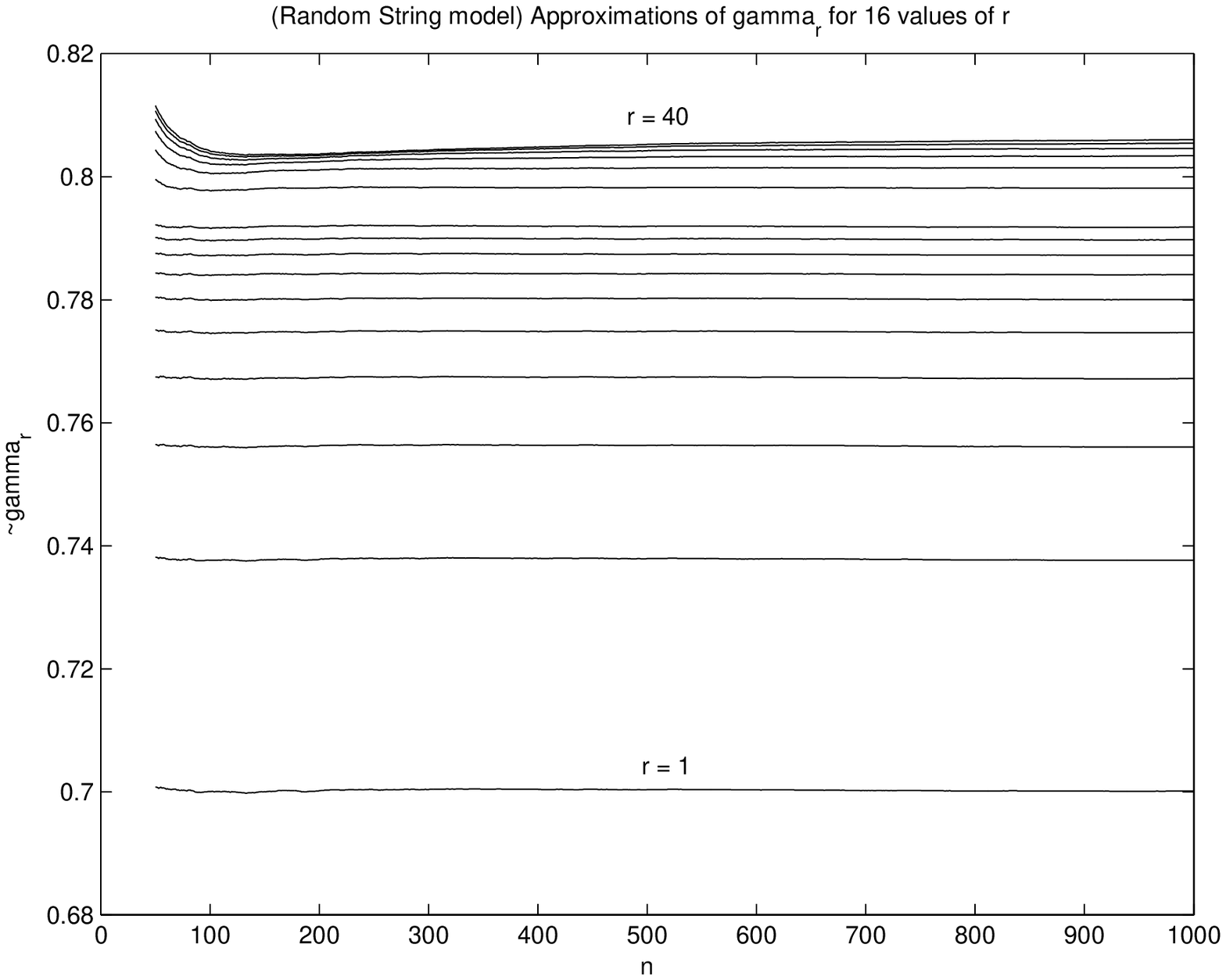}
\end{center}
\caption{The Monte Carlo approximations of $\frac{\mathbf{EL}_{n,2,r}}{n}$
and this same data corrected by (\protect\ref{finitesizefix})\ to obtain the
limiting behavior. }
\label{twographs}
\end{figure}

For the Bernoulli Matching model we also can compute $\mathbf{EL}_{n,2,r}^{B}
$ exactly for small $n$ by applying (\ref{recur}) directly beginning with $%
\overrightarrow{P_{r}(z)}$. This allows us to do two checks on the quality
of the Monte Carlo approximations. Firstly, we can observe the difference $%
\frac{MonteCarlo(\mathbf{EL}_{n,2,r}^{B})}{n}-\frac{\mathbf{EL}_{n,2,r}^{B}}{%
n}$. The statistic
\begin{equation*}
S_{r}\equiv \frac{1}{1000}\dsum_{j=1}^{1000}\left( \frac{MonteCarlo(\mathbf{%
EL}_{n,2,r}^{B})}{n}-\frac{\mathbf{EL}_{n,2,r}^{B}}{n}\right) ^{2}
\end{equation*}%
gives us an idea of how crude an approximation we get with $10000$ trials.
Also, we can see how good the approximation $\mathbf{EL}_{n,2,r}^{B}$ $%
\mathbf{\sim }$ $\gamma _{2,r}^{B}n-A_{r}^{B}$ is by using that on the exact
values of $\mathbf{EL}_{n,2,r}^{B}$ to extrapolate $\gamma _{2,r}^{B}$ (for
this extrapolation we use $n=1...2000$).%
\begin{eqnarray*}
&&%
\begin{tabular}{l|llll}
$r$ & $S_{r}$ & Monte Carlo $\gamma _{2,r}^{B}$ & $\gamma _{2,r}^{B}$ from $%
\mathbf{EL}_{n,2,r}^{B}$ & $\gamma _{2,r}^{B}$ $\QATOP{\text{from fractions}%
}{\text{derived previously}}$ \\ \hline
$1$ & $5.2994\times 10^{-8}$ & $0.7272634$ & $0.7272727273$ & $0.7272727272$
\\ \hline
$2$ & $5.0758\times 10^{-8}$ & $0.7716676$ & $0.7715736043$ & $0.7715736040$
\\ \hline
$3$ & $1.0180\times 10^{-8}$ & $0.7889874$ & $0.7889693851$ & $0.7889693853$
\\ \hline
$4$ & $1.5954\times 10^{-8}$ & $0.7981354$ & $0.7982222051$ & $-$%
\end{tabular}
\\
&&%
\begin{tabular}{l|lll}
$r$ & Monte Carlo $A_{r}^{B}$ & $A_{r}^{B}$ from $\mathbf{EL}_{n,2,r}^{B}$ &
$A_{r}^{B}\QATOP{\text{from fractions}}{\text{derived previously}}$ \\ \hline
$1$ & $0.2771$ & $0.264463$ & $0.2644628$ \\ \hline
$2$ & $0.4626$ & $0.434745$ & $0.4347445$ \\ \hline
$3$ & $0.5641$ & $0.574312$ & $-$ \\ \hline
$4$ & $0.6852$ & $0.696534$ & $-$%
\end{tabular}%
\end{eqnarray*}

We also note that $MonteCarlo(\gamma _{2,1})=0.7001417$ compared to $\gamma
_{2,1}=.7$ and $MonteCarlo(A_{2,1})=0.2652$ compared to $A_{2,1}=.28$

\section{Conclusions and future work}

It is hoped that the results presented in this paper lead the way to more
significant results. In particular, it is hoped that the Random String model
analysis may lead to a short proof of $\gamma _{2,1}<\gamma _{2,1}^{B}$. The
limiting behavior of $\gamma _{k,1}-\gamma _{k,1}^{B}$ would also be of
interest. We seek a conjecture for the quantities $\gamma _{k,r}^{B}$,
though it is unclear if trying to determine $\gamma _{k}^{B}$ via $%
\lim_{r\longrightarrow \infty }\gamma _{k,r}^{B}=\gamma _{k}^{B}$ is a good
idea.

The pseudoproof of $\gamma _{k}^{B}=\frac{2}{1+\sqrt{k}}$ given by Boutet de
Monvel may provide a way to simplify the r-reach computations. The limiting
behavior of r-reach may be describable only by differences between adjacent
values of $\mathbf{R,}$ thereby reducing the "problemsize" from $2^{2r}$ to $%
2r$. Preliminary investigations suggest that this reduction may be possible
but not as straight forward as the argument in the pseudoproof.

\section{Appendix}

The expanded version of the Bernoulli Matching model $r=1,k=2$ case has
matrices as follows. These are given for comparison with the matrices for
the Random String model $r=1,k=2$ case.

\begin{equation*}
\text{ \ \ \ \ }%
\begin{array}{l}
{\small P}_{n-1}^{\text{off}}{\small (z,z,z)} \\
{\small P}_{n-1}^{\text{off}}{\small (z,z,z-1)} \\
{\small P}_{n-1}^{\text{off}}{\small (z,z-1,z)} \\
{\small P}_{n-1}^{\text{off}}{\small (z,z-1,z-1)} \\
{\small P}_{n-1}^{\text{on}}{\small (z,z,z)} \\
{\small P}_{n-1}^{\text{on}}{\small (z,z,z-1)} \\
{\small P}_{n-1}^{\text{on}}{\small (z,z-1,z)} \\
{\small P}_{n-1}^{\text{on}}{\small (z,z-1,z-1)}%
\end{array}%
\begin{bmatrix}
{\small 1/8} & {\small 0} & {\small 0} & {\small 0} & {\small 0} & {\small 0}
& {\small 0} & {\small 0} \\
{\small 1/4} & {\small 0} & {\small 0} & {\small 0} & {\small 0} & {\small 0}
& {\small 0} & {\small 0} \\
{\small 1/4} & {\small 0} & {\small 0} & {\small 0} & {\small 0} & {\small 0}
& {\small 0} & {\small 0} \\
{\small 1/2} & {\small 0} & {\small 0} & {\small 0} & {\small 0} & {\small 0}
& {\small 0} & {\small 0} \\
{\small 1/8} & {\small 0} & {\small 0} & {\small 0} & {\small 0} & {\small 0}
& {\small 0} & {\small 0} \\
{\small 1/4} & {\small 0} & {\small 0} & {\small 0} & {\small 0} & {\small 0}
& {\small 0} & {\small 0} \\
{\small 1/4} & {\small 0} & {\small 0} & {\small 0} & {\small 0} & {\small 0}
& {\small 0} & {\small 0} \\
{\small 1/2} & {\small 0} & {\small 0} & {\small 0} & {\small 0} & {\small 0}
& {\small 0} & {\small 0}%
\end{bmatrix}%
=\mathbf{M}
\end{equation*}%
\begin{equation*}
\text{ \ \ \ \ }%
\begin{array}{l}
{\small P}_{n-1}^{\text{off}}{\small (z-1,z-1,z-1)} \\
{\small P}_{n-1}^{\text{off}}{\small (z-1,z-1,z-2)} \\
{\small P}_{n-1}^{\text{off}}{\small (z-1,z-2,z-1)} \\
{\small P}_{n-1}^{\text{off}}{\small (z-1,z-2,z-2)} \\
{\small P}_{n-1}^{\text{on}}{\small (z-1,z-1,z-1)} \\
{\small P}_{n-1}^{\text{on}}{\small (z-1,z-1,z-2)} \\
{\small P}_{n-1}^{\text{on}}{\small (z-1,z-2,z-1)} \\
{\small P}_{n-1}^{\text{on}}{\small (z-1,z-2,z-2)}%
\end{array}%
\begin{bmatrix}
{\small 1/8} & {\small 1/8} & {\small 1/8} & {\small 0} & {\small 1/8} &
{\small 1/8} & {\small 1/8} & {\small 1/8} \\
{\small 0} & {\small 1/4} & {\small 0} & {\small 0} & {\small 0} & {\small %
1/4} & {\small 0} & {\small 1/4} \\
{\small 0} & {\small 0} & {\small 1/4} & {\small 0} & {\small 0} & {\small 0}
& 1{\small /4} & {\small 1/4} \\
{\small 0} & {\small 0} & {\small 0} & {\small 0} & {\small 0} & {\small 0}
& {\small 0} & {\small 1/2} \\
{\small 1/8} & {\small 1/8} & {\small 1/8} & {\small 0} & {\small 1/8} &
{\small 1/8} & {\small 1/8} & {\small 1/8} \\
{\small 0} & 1/4 & {\small 0} & {\small 0} & {\small 0} & 1/4 & {\small 0} &
{\small 1/4} \\
{\small 0} & {\small 0} & {\small 1/4} & {\small 0} & {\small 0} & {\small 0}
& {\small 1/4} & {\small 1/4} \\
{\small 0} & {\small 0} & {\small 0} & {\small 0} & {\small 0} & {\small 0}
& {\small 0} & 1/2%
\end{bmatrix}%
=\mathbf{N}
\end{equation*}

\pagebreak

\end{document}